\documentclass[english]{article}
\usepackage[T1]{fontenc}
\usepackage[latin9]{inputenc}
\usepackage{geometry}
\usepackage{amssymb}
\usepackage{amsmath}
\usepackage{graphicx}

\makeatletter
\usepackage{dsfont}

\makeatother

\usepackage[english,francais]{babel}
\addto\extrasfrench{%
   \providecommand{\fg}{\ifdim\lastskip>\z@\unskip\fi~\frqq}%
}

\newtheorem{theorem}{Theorem}[section]
\newtheorem{lemma}[theorem]{Lemma}
\newtheorem{e-proposition}[theorem]{Proposition}

\newtheorem{e-definition}[theorem]{Definition\rm}

\newtheorem{example}{\it Example\/}

\begin{document}
\selectlanguage{english}
\title{Large deviations of a velocity jump process with a Hamilton-Jacobi approach}

\author{Nils Caillerie}\maketitle

\begin{abstract}
\selectlanguage{english}
% Text of abstract in English
We study a random process on $\mathbb{R}^{n}$ moving in straight
lines and changing randomly its velocity at random exponential times.
We focus more precisely on the Kolmogorov equation in the hyperbolic
scale $\left(t,x,v\right)\to\left(\frac{t}{\varepsilon},\frac{x}{\varepsilon},v\right)$,
with $\varepsilon>0$, before proceeding to a Hopf-Cole transform,
which gives a kinetic equation on a potential. We show convergence
as $\varepsilon\to0$ of the potential towards the viscosity solution
of a Hamilton-Jacobi equation $\partial_{t}\varphi+H\left(\nabla_{x}\varphi\right)=0$
where the hamiltonian may lack $\mathcal{C}^{1}$ regularity, which is
quite unseen in this type of studies.

\vskip 0.5\baselineskip

\selectlanguage{francais}
% Text of abstract in French
\begin{centering}\noindent{\bf R\'esum\'e} \vskip 0.5\baselineskip \noindent \end{centering}
\textbf{Grandes d\'eviations pour un processus \`a sauts de vitesse avec une approche Hamilton-Jacobi}. Nous nous int\'eressons \`a un processus al\'eatoire sur $\mathbb{R}^{n}$
qui alterne des phases de mouvements rectilignes uniformes et change
de vitesse \`a des temps exponentiels. Nous \'etudions plus pr\'ecis\'ement
l'\'equation de Kolmogorov apr\`es r\'e\'echelonnement hyperbolique $\left(t,x,v\right)\to\left(\frac{t}{\varepsilon},\frac{x}{\varepsilon},v\right)$,
$\varepsilon>0$, puis nous effectuons une transform\'ee de Hopf-Cole
qui nous donne une \'equation cin\'etique suivie par un potentiel. Nous
montrons la convergence pour $\varepsilon\to0$ de ce potentiel vers
la solution de viscosit\'es d'une \'equation de Hamilton-Jacobi
$\partial_{t}\varphi+H\left(\nabla_{x}\varphi\right)=0$ o\`u le hamiltonien peut
pr\'esenter une singularit\'e $\mathcal{C}^{1}$, ce qui est assez in\'edit
dans ce type d'\'etudes.

\end{abstract}

% now the Version fran\'eaise abr\'eg\'ee, if it exists
\selectlanguage{francais}
\section*{Version fran\c{c}aise abr\'eg\'ee}
% Text of your Version fran\'eaise abr\'eg\'ee here.
% Note you do not need to repeat here equations that you use in the
% main text - for example 'voir (3)' is quite acceptable.

Nous nous donnons une densit\'e de probabilit\'e $M\in L^{1}\left(\mathbb{R}^{n}\right)$
et nous notons $V$ son support. Nous supposons que $V$ est compact
et que $0$ appartient \`a l'int\'erieur de l'enveloppe convexe de $V$,
que l'on note $\mathrm{Conv}\left(V\right)$. Pour $p\in\mathbb{R}^{n}$,
nous notons $\mu\left(p\right)=\mathrm{max}\left\{ v\cdot p\mid v\in\mathrm{Conv}\left(V\right)\right\} $.
Nous \'etudions le mouvement de particules dans $\mathbb{R}^{n}$ suivant
le processus de Markov d\'eterministe par morceaux d\'efini comme suit
: une particule donn\'ee se d\'eplace de mani\`ere rectiligne uniforme avec
une vitesse $v\in V$ tir\'ee al\'eatoirement en suivant la loi de probabilit\'e
$M\left(v'\right)dv'$. \`A des temps exponentiels de param\`etre 1, la
particule change de direction en tirant une nouvelle vitesse tir\'ee
selon la loi $M\left(v'\right)dv'$. Afin d'\'etudier des r\'esultats de larges d\'eviations du processus similairement aux techniques d\'evelopp\'ees dans \cite{Bressloff}-\cite{Faggionato}, nous nous int\'eressons \`a l'\'equation
de Chapman-Kolmogorov forward suivie par la densit\'e de particules apr\`es un r\'e\'echelonement
hyperbolique $\left(t,x,v\right)\to\left(\frac{t}{\varepsilon},\frac{x}{\varepsilon},v\right)$,
$\varepsilon>0$ : 
\[
\partial_{t}f^{\varepsilon}+v\cdot\nabla_{x}f^{\varepsilon}=\frac{1}{\varepsilon}\left(M\left(v\right)\rho^{\varepsilon}-f^{\varepsilon}\right),\quad\left(t,x,v\right)\in\mathbb{R}_{+}\times\mathbb{R}^{n}\times V.
\]
 Nous \'etudions plus particuli\`erement l'\'equation v\'erifi\'ee par un potentiel
$\varphi^{\varepsilon}$ obtenu apr\`es passage par une transform\'ee
de Hopf-Cole : $f^{\varepsilon}\left(t,x,v\right)=M\left(v\right)e^{-\frac{\varphi^{\varepsilon}\left(t,x,v\right)}{\varepsilon}}$.
Nous cherchons alors une \'eventuelle limite pour $\varphi^{\varepsilon}$.
Nous proc\'edons \`a un d\'eveloppement WKB : $\varphi^{\varepsilon}=\varphi+\varepsilon\eta$,
ce qui am\`ene, en posant $p=\nabla_{x}\varphi$ et $H=-\partial_{t}\varphi$,
\`a la r\'esolution d'un probl\`eme spectral dans l'espace des mesures positives
: chercher $\left(H,Q\right)$ un couple valeur/vecteur propres associ\'e
\`a l'op\'erateur $Q\mapsto\left(v\cdot p-1\right)Q+\int_{V} M'Q'dv'$. On obtient une \'equation de Hamilton-Jacobi $\partial_{t}\varphi+H\left(\nabla_{x}\varphi\right)=0$.
Pour $n=1$ et $M\geq\delta>0$ sur son support, le vecteur propre $Q$ a une densit\'e
et conduit \`a un hamiltonien $H$ d\'efini par l'\'equation implicite
\[
\int_{V}\frac{M\left(v\right)}{1+H\left(p\right)-v\cdot p}dv=1.
\]
La positivit\'e de $Q$ garantit que $H\left(p\right)\geq\mu\left(p\right)-1$. En dimension sup\'erieure toutefois, et m\^eme si $M\geq\delta>0$, cette
\'equation peut ne pas avoir de solution $H\left(p\right)$ lorsque
$p$ devient grand. Cela se manifeste pour le vecteur propre par une
concentration de la mesure $Q$ autour des valeurs $v$ qui annulent $1+H\left(p\right)-v\cdot p$,
ce qui force $H\left(p\right)=\mu\left(p\right)-1$. Cette transition entra\^ine une singularit\'e $\mathcal{C}^1$ du hamiltonien.

Nous d\'emontrons la convergence de $\varphi^{\varepsilon}$ vers $\varphi$,
o\`u $\varphi$ est solution de viscosit\'e \cite{Crandall} de
l'\'equation de Hamilton-Jacobi en utilisant la m\'ethode de la fonction test perturb\'ee \cite{perturbed}.

\selectlanguage{english}
% main text
\section{Introduction}
\label{introduction}
% etc, etc

We continue the work initiated in \cite{Bouin}-\cite{BC12}. Let $M\in L^{1}\left(\mathbb{R}^{n}\right)$
be a probability density function. We suppose that the support
of $M$, which we denote $V$, is compact and that $0$ belongs to
the interior $ $of $\mathrm{Conv}\left(V\right)$, the convex hull
of $V$. We denote by $\left|\cdot\right|$ the euclidian
norm in $\mathbb{R}^{n}$ and by $\cdot$ the canonical scalar product.
For $p\in\mathbb{R}^{n}$, we define
\begin{equation}
\mu\left(p\right):=\mathrm{max}\left\{ v\cdot p\mid v\in\mathrm{Conv}\left(V\right)\right\},\label{eq:mu}
\end{equation}
$\mathrm{Arg}\mu\left(p\right):=\left\{ v\in\mathrm{Conv}\left(V\right)\mid v\cdot p=\mu\left(p\right)\right\} $ and $\mathrm{Sing}\left(M\right):={\left\{ p\in\mathbb{R}^{n},\quad{\int_{V}}\frac{M\left(v\right)}{\mu\left(p\right)-v\cdot p}dv\leq1\right\}}$.

We focus on the motion dynamics in $\mathbb{R}^{n}$ of particles given by
the following piecewise deterministic Markov process: a particle moves successively
in straight lines with velocity $v$, chosen randomly with probability
distribution $M\left(v'\right)dv'$. At random exponential times (with
parameter $1$), the particle changes its velocity, choosing randomly
a new velocity with distribution $M\left(v'\right)dv'$. The Chapman-Kolmogorov forward
equation associated to the probability density
function $f\left(t,x,v\right)$ of this process is given by: 
\begin{equation}
\partial_{t}f+v\cdot\nabla_{x}f=M\rho-f,\quad\left(t,x,v\right)\in\mathbb{R}_{+}\times\mathbb{R}^{n}\times V,
\label{eq:1}
\end{equation}
where $\rho\left(t,x\right)=\int_{V}f\left(t,x,v\right)dv$. In order to investigate large deviation principles for the process, one can study the large scale hyperbolic limit $\left(t,x\right)\to\left(\frac{t}{\varepsilon},\frac{x}{\varepsilon}\right)$
with $\varepsilon>0$. In this scale, the kinetic equation (\ref{eq:1})
reads:
\begin{equation}
\partial_{t}f^{\varepsilon}+v\cdot\nabla_{x}f^{\varepsilon}=\frac{1}{\varepsilon}\left(M\rho^{\varepsilon}-f^{\varepsilon}\right),\quad\left(t,x,v\right)\in\mathbb{R}_{+}\times\mathbb{R}^{n}\times V.\label{eq:scaling}
\end{equation}
Then, we perform the following Hopf-Cole transformation: $f^{\varepsilon}\left(t,x,v\right)=M\left(v\right)e^{-\frac{\varphi^{\varepsilon}\left(t,x,v\right)}{\varepsilon}}$,
where we expect the potential $\varphi^{\varepsilon}$ to become independent
of $v$ as $\varepsilon\to0$. Such techniques have already been studied for a more general case of Markov process with a finite discrete set of states in \cite{Bressloff} and, from a probabilistic point of view, in \cite{Faggionato}. Here, assume that the initial condition
is well-prepared, i.e. it does not depend on $v$ and $\varepsilon$: $\varphi^{\varepsilon}\left(0,x,v\right)=\varphi_{0}\left(x\right)$. The equation satisfied
by $\varphi^{\varepsilon}$ reads
\begin{equation}
\partial_{t}\varphi^{\varepsilon}+v\cdot\nabla_{x}\varphi^{\varepsilon}=\int_{V}M\left(v'\right)\left(1-e^{\frac{\varphi^{\varepsilon}-\varphi'^{\varepsilon}}{\varepsilon}}\right)dv',\quad\left(t,x,v\right)\in\mathbb{R}_{+}\times\mathbb{R}^{n}\times V.
\label{eq:kinetic}
\end{equation}
As in \cite{Wentzell}, the limit potential satisfy a Hamilton-Jacobi equation. Surprisingly enough, our Hamiltonian may lack $\mathcal{C}^{1}$ regularity as we will show in Proposition \ref{prop:Hamiltonien}.
\begin{theorem}
\label{thm:Let--be}Under the previous assumptions, $\varphi^{\varepsilon}$
converges locally uniformly on $\mathbb{R}_{+}\times\mathbb{R}^{n}\times V$
toward $\varphi$, where $\varphi$ does not depend on $v$. Moreover,
$\varphi$ is the viscosity solution of the following Hamilton-Jacobi
equation:
\begin{equation}
\partial_{t}\varphi\left(t,x\right)+H\left(\nabla_{x}\varphi\left(t,x\right)\right)=0,\quad\left(t,x\right)\in\mathbb{R}_{+}\times\mathbb{R}^{n},\label{eq:hamilton Jacobi}
\end{equation}
where the hamiltonian $H$ is given as follows: if $p\in\mathrm{Sing}\left(M\right)$,
then $H\left(p\right)=\mu\left(p\right)-1$. Else, $H\left(p\right)$ is uniquely determined by the following formula:
\begin{equation}
\int_{V}\frac{M\left(v\right)}{1+H\left(p\right)-v\cdot p}dv=1.\label{eq:hamiltonienimplicite}
\end{equation}\end{theorem}

\section{Identification of the hamiltonian\label{sec:2}}

In order to identify the limit $\varepsilon\to0$ of the equation
(\ref{eq:kinetic}) we perform the formal WKB expansion: $\varphi^{\varepsilon}\left(t,x,v\right)=\varphi\left(t,x\right)+\varepsilon\eta\left(t,x,v\right),$
where $\varphi$ and $\eta$ are to be determined. Plugging this ansatz
into the kinetic formulation (\ref{eq:kinetic}), we get, taking the
formal limit $\varepsilon\to0$:
\[
\partial_{t}\varphi+v\cdot\nabla_{x}\varphi=1-\int_{V}M\left(v'\right)e^{\eta-\eta'}dv'.
\]
Let us write $p=\nabla_{x}\varphi$ and $H=-\partial_{t}\varphi$.
The equation for $Q=e^{-\eta}$ is the following spectral problem: $HQ=\left(v\cdot p-1\right)Q+\int_{V}M\left(v'\right)Q\left(v'\right)dv'\label{eq:pbspectral}$.
The positivity of $Q$ yields $H\geq v\cdot p-1$ for
all $v\in V$ hence $H\geq\mu\left(p\right)-1$. Suppose
$H>\mu\left(p\right)-1$. Then, $1+H-vp>0$
for all $v\in V$ and ${\displaystyle Q\left(v\right)=\frac{\int_{V}M\left(v'\right)Q\left(v'\right)dv'}{1+H-v\cdot p}}$.
Integrating against $M$ with respect to $v$, we obtain the following
problem: find $H$ such that ${\int_{V}\frac{M\left(v\right)}{1+H-v\cdot p}dv=1}$.
If $p\in\mathrm{Sing}\left(M\right)^{c}$, by monotonicity, such $H$
exists and is unique. Equation (\ref{eq:hamiltonienimplicite}), however,
does not have a solution for $p\in\mathrm{Sing}\left(M\right)$, so
we necessarily have $H=\mu\left(p\right)-1$. Then,
a possible solution of the spectral problem is the positive measure
$Q=\frac{dv}{\mu\left(p\right)-v\cdot p}+\alpha\left(p\right)\delta_{w}$
where $\alpha\left(p\right)=1-{\int_{V}\frac{M\left(v\right)}{\mu\left(p\right)-vp}dv}\geq0$
and $\delta_{w}$ is the Dirac measure centered in
$w\in\mathrm{Arg}\mu\left(p\right)\cap V$.

Here is an example where $\mathrm{Sing}\left(M\right)\neq\emptyset$:
\begin{example}
\label{(partie-singuliere)} Let $n>1$ and $M=\omega_{n}^{-1}.\mathds{1}_{\overline{B\left(0,1\right)}}$
where $\omega_{n}$ is the Lebesgue measure of the $n$-dimensional
unit ball. Then, $\mathrm{Sing}\left(M\right)=B\left(0,\frac{n}{n-1}\right)^{c}$. Indeed,
for $p=\left|p\right|\cdot e_{1}$, we have $\mu\left(p\right)=\left|p\right|$
and $v\cdot p=\left|p\right|v_{1}$ hence
\[
\int_{V}\frac{M\left(v\right)}{\mu\left(p\right)-v\cdot p}dv=\frac{1}{\left|p\right|\omega_{n}}\int_{B\left(0,1\right)}\frac{1}{1-v_{1}}dv=\frac{\omega_{n-1}}{\left|p\right|\omega_{n}}\int_{-1}^{1}\frac{\left(1-v_{1}^{2}\right)^{\frac{n-1}{2}}}{1-v_{1}}dv_{1}=\frac{1}{\left|p\right|}\times\frac{n}{n-1}.
\]
By rotational invariance, we conclude that $\mathrm{Sing}\left(M\right)=B\left(0,\frac{n}{n-1}\right)^{c}$.
The Figure 1 gives illustrations of the hamiltonian and $\mu$ as functions of the radius of $p$, 
in the cases $n=1$ and $n=3$. In the cases $n=3$ we can see the $\mathcal{C}^{1}$ singularity where $\left|p\right|=\frac{3}{2}$.
\begin{figure}[h]
\caption{Blue plain lines : Hamiltonian for $n=1,3$ and $M=\omega_{n}^{-1}.\mathds{1}_{\overline{B\left(0,1\right)}}$.
Black dotted lines : $\left|p\right|\mapsto\mu\left(p\right)-1$.
Lignes pleines bleues : Hamiltonien pour $n=1,3$ et $M=\omega_{n}^{-1}.\mathds{1}_{\overline{B\left(0,1\right)}}$.
Lignes noires en pointill\'es : $\left|p\right|\mapsto\mu\left(p\right)-1$.}

\begin{centering}
\includegraphics[width=0.48\linewidth=0.48]{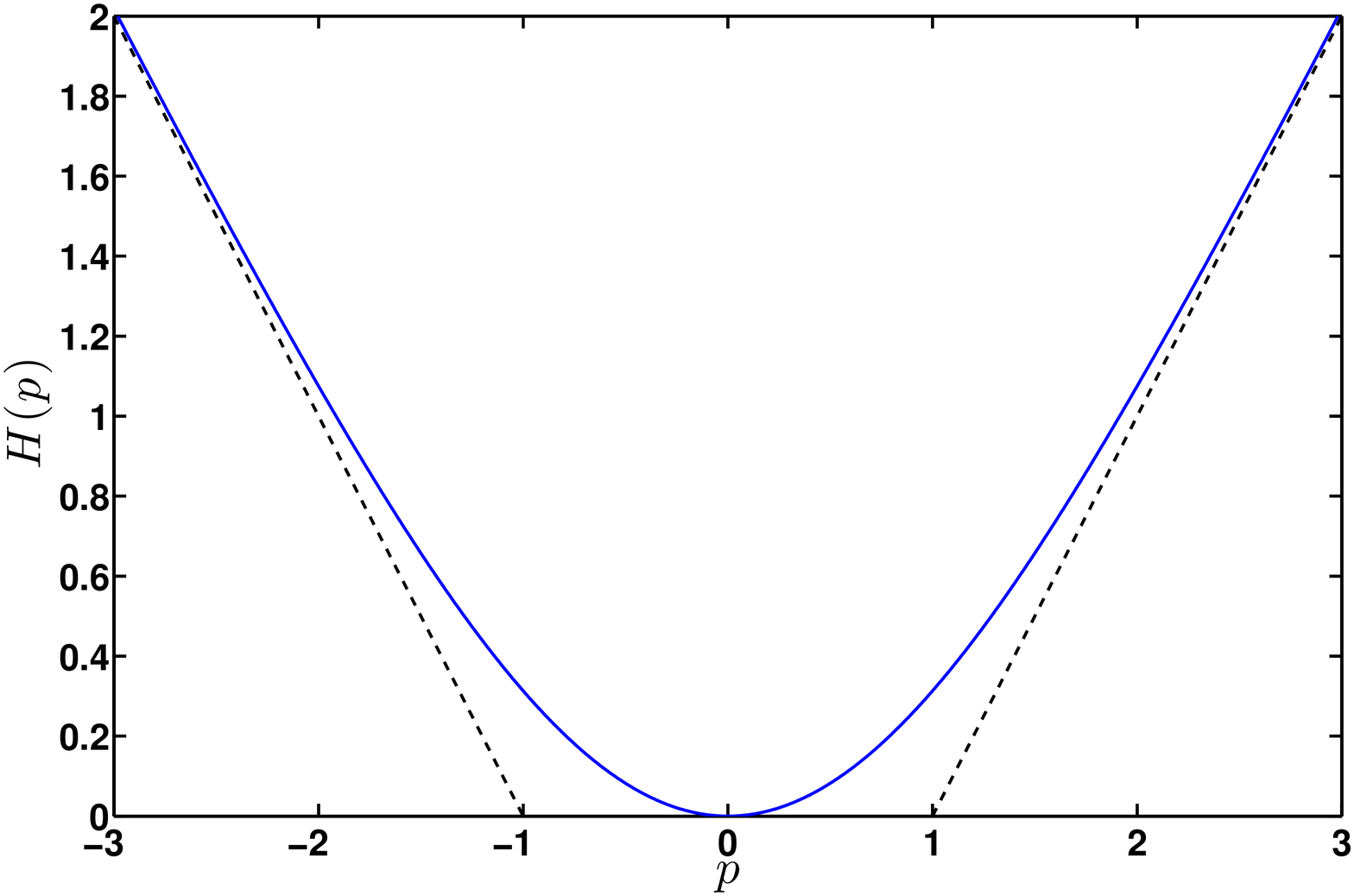}\includegraphics[width=0.48\linewidth=0.48]{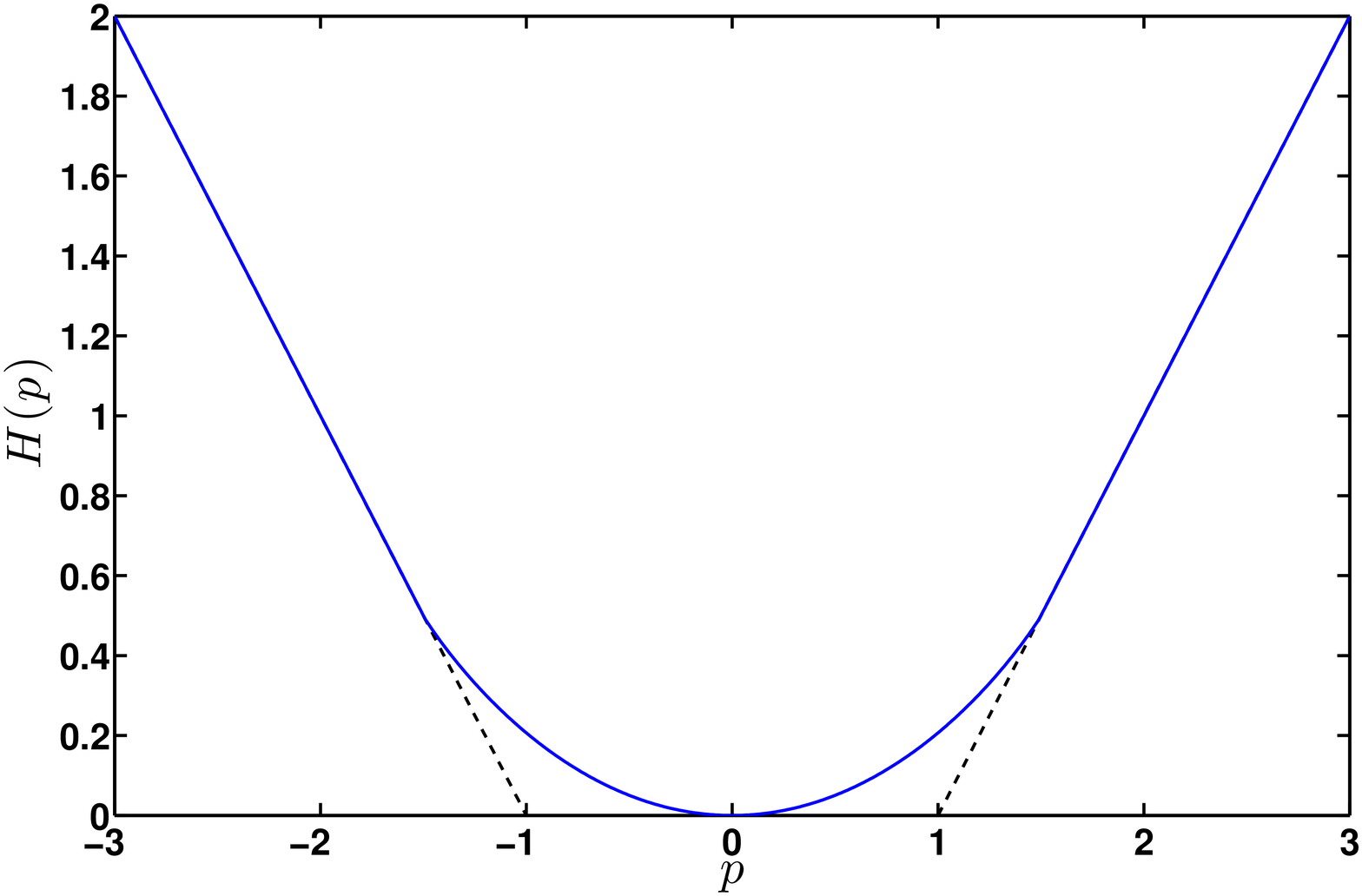}
\par\end{centering}

\end{figure}
\end{example}
\begin{e-proposition}
\label{prop:Hamiltonien}The following properties hold:

(i) The set $\mathrm{Sing}\left(M\right)^{c}$ is convex.

(ii) The function $H$ is continuous and convex. 

(iii) If $\mathrm{Sing}\left(M\right)\neq\emptyset$, then $H$ is
not $\mathcal{C}^{1}$. More precisely, $\nabla H$ has a jump discontinuity at $\partial\mathrm{Sing}M$.\end{e-proposition}

\textbf{Proof} Let us first notice that $\mu$ is positively $1$-homogeneous. Moreover, it is convex since it is a
supremum of linear functions.

(i) Let $p,q\in\mathrm{Sing}\left(M\right)^{c}$ with $p\neq q$.
Since $\mu$ is convex, we have for all $\tau\in\left[0,1\right]$
\begin{eqnarray*}
I\left(\tau\right):=\int_{V}\frac{M\left(v\right)}{\mu\left(p\right)-v\cdot p+\tau\left(\mu\left(q\right)-\mu\left(p\right)-v\cdot\left(q-p\right)\right)}dv & \leq & \int_{V}\frac{M\left(v\right)}{\mu\left(\left(1-\tau\right)p+\tau q\right)-v\cdot\left(\left(1-\tau\right)p+\tau q\right)}dv.
\end{eqnarray*}
Moreover, $I\left(0\right),I\left(1\right)>1$ and $I$ is differentiable on $\left[0,1\right]$ with
\begin{eqnarray*}
\partial_{\tau}I\left(\tau\right) & = & \int_{V}\frac{M\left(v\right)}{\left(\mu\left(p\right)-v\cdot p+\tau\left(\mu\left(q\right)-\mu\left(p\right)-v\cdot\left(q-p\right)\right)\right)^{2}}\left(\mu\left(p\right)-\mu\left(q\right)-v\cdot\left(p-q\right)\right)dv.
\end{eqnarray*}
It is clear that the sign of $\partial_{\tau}I$ does not change
hence $I \left(\tau\right)>1$, which proves (i).

(ii) We refer to \cite{BC12} to prove that $H$ is $\mathcal{C}^{2}$
and strictly convex on $\mathrm{Sing}\left(M\right)^{c}$ and that
\begin{equation}
{\displaystyle \int_{V}\frac{M\left(v\right)}{\left(1+H\left(q\right)-v\cdot q\right)^{2}}\left(\nabla H\left(q\right)-v\right)dv}=0,\quad\forall q\in\mathrm{Sing}\left(M\right)^{c}.\label{eq:gradH}
\end{equation}
In particular, $\nabla H\left(q\right)\in\mathrm{Conv}\left(V\right)$
for all $q\in\mathrm{Sing}\left(M\right)^{c}$. It is easy to see
that $H$ is continuous in the interior of $\mathrm{Sing}\left(M\right)$. To show continuity of $H$ on $\partial\mathrm{Sing}\left(M\right)$,
let $\left(p_{m}\right)_{m}$ converge to $p\in\partial\mathrm{Sing}\left(M\right)\subset\mathrm{Sing}\left(M\right)$.
If we can extract a subsequence $\left(p_{m_{l}}\right)_{l}\subset\mathrm{Sing}\left(M\right)$,
then $H\left(p_{m_{l}}\right)=\mu\left(p_{m_{l}}\right)-1\underset{l\to\infty}{\longrightarrow}\mu\left(p\right)-1=H\left(p\right)$.
If not, then $p_{m}\in\mathrm{Sing}\left(M\right)^{c}$ for $m$ large enough and $1=\int_{V}\frac{M\left(v\right)}{1+H\left(p_{m}\right)-v\cdot p_{m}}dv<\int_{V}\frac{M\left(v\right)}{\mu\left(p_{m}\right)-v\cdot p_{m}}dv$.
Taking the limit, we get by dominated convergence $1=\int_{V}\underset{m\to\infty}{\mathrm{lim}}\frac{M\left(v\right)}{1+H\left(p_{m}\right)-v\cdot p_{m}}dv\leq\int_{V}\frac{M\left(v\right)}{\mu\left(p\right)-v\cdot p_{m}}dv\leq1$ hence $H\left(p_{m}\right)\underset{m\to\infty}{\longrightarrow}\mu\left(p\right)-1=H\left(p\right)$.

We now show that $H$ is convex by proving that it is a maximum of
convex functions:
\begin{eqnarray}
H\left(p\right) & = & \mathrm{max}\left(\mathrm{sup}\left\{ \nabla H\left(q\right)\cdot\left(p-q\right)+H\left(q\right)\mid q\in\mathrm{Sing}\left(M\right)^{c}\right\} ,\mu\left(p\right)-1\right),\quad\forall p\in\mathbb{R}^{n}.\label{eq:supH}
\end{eqnarray}
In $\mathrm{Sing}\left(M\right)^{c}$, (\ref{eq:supH}) holds
by convexity of $H$ and $H\left(p\right)>\mu\left(p\right)-1$.
Let $p\in\mathrm{Sing}\left(M\right)$ and $q\in\mathrm{Sing}\left(M\right)^{c}$.
By convexity of $\mathrm{Sing}\left(M\right)^{c}\ni0$, there exists a
unique $\lambda\in\left(0,1\right]$ such that $\lambda p\in\partial\mathrm{Sing}\left(M\right)$.
For all $\tau\in\left[0,1\right]$, we set $\omega_{1}\left(\tau\right):=\mu\left(\tau p\right)-1=\tau\mu\left(p\right)-1$ and $\omega_{2}\left(\tau\right):=\nabla H\left(q\right)\cdot\left(\tau p-q\right)+H\left(q\right).$ By continuity of $H$, $\mu\left(\lambda p\right)-1=H\left(\lambda p\right)\geq\nabla H\left(q\right)\left(\lambda p-q\right)+H\left(q\right)$
hence $\omega_{1}\left(\lambda\right)\geq\omega_{2}\left(\lambda\right)$.
Moreover, $\omega_{1}$ and $\omega_{2}$ are both differentiable
and $\partial_{\tau}\omega_{1}\left(\tau\right)=\mu\left(p\right)\geq\nabla H\left(q\right)\cdot p=\partial_{\tau}\omega_{2}\left(\tau\right)$
since $\nabla H\left(q\right)\in\mathrm{Conv}\left(V\right)$. Hence,
$\omega_{1}\left(1\right)\geq\omega_{2}\left(1\right)$, which ends
the proof of (ii).

(iii) Suppose $\mathrm{Sing}\left(M\right)\neq\emptyset$ and $H$
is $\mathcal{C}^{1}$. Since $H+1=\mu$ is positive homogeneous of
degree 1 on $\mathrm{Sing}\left(M\right)$ and since $\lambda p\in\mathrm{Sing}\left(M\right)$ for all $\lambda\geq1$ and $p\in\mathrm{Sing}\left(M\right)$, we know that $\nabla H\left(p\right)\cdot p=H\left(p\right)+1=\mu\left(p\right)$
for all $p\in\mathrm{Sing}\left(M\right)\subset\mathrm{Sing}\left(M\right)$
hence $p\cdot\left(\nabla H\left(p\right)-v\right)\geq0$, for all $v\in V$, the inequality being strict on a neighborhood of $0$. Then,
\begin{equation}
p\cdot{\displaystyle \int_{V}\frac{M\left(v\right)}{\left(1+H\left(p\right)-v\cdot p\right)^{2}}\left(\nabla H\left(p\right)-v\right)dv}>0,\quad\forall p\in\partial\mathrm{Sing}\left(M\right).\label{eq:contra}
\end{equation}
By continuity, equations (\ref{eq:gradH}) and (\ref{eq:contra})
are contradictory.\begin{flushright}$\square$\end{flushright}

\section{Proof of Theorem \ref{thm:Let--be}}

Let $\varphi_{0}\in W^{1,\infty}\left(\mathbb{R}^{n}\right)$. We refer to Proposition 2.1 in \cite{BC12} to prove that the Cauchy Problem (\ref{eq:kinetic}) with initial condition $\varphi_{0}$
has a unique solution $\varphi^{\varepsilon}\in W^{1,\infty}$ which is locally
(in $t$) uniformly (in $\varepsilon$, $x$ and $v$) bounded in
norm $W^{1,\infty}$. In particular, let us mention that
\begin{equation}
0\leq\varphi^{\varepsilon}\left(t,\cdot,\cdot\right)\leq\left\Vert \varphi_{0}\right\Vert _{\infty}\label{apriori:Linfty},\quad
\left\Vert \nabla_{v}\varphi^{\varepsilon}\left(t,\cdot,\cdot\right)\right \Vert _{\infty}\leq\left\Vert \nabla_{x}\varphi_{0}\right\Vert _{\infty}.
\end{equation}
Using the Arzel\'a-Ascoli theorem, we extract a locally uniformly converging
subsequence. We denote by $\varphi$ the limit. The function $\varphi$
does not depend on $v$ since ${\int_{V}}M\left(v\right)e^{\frac{\varphi^{\varepsilon}-\varphi'^{\varepsilon}}{\varepsilon}}dv$
is uniformly bounded on $\left[0,T\right]\times\mathbb{R}^{n}\times V$
for all $T>0$. We use the perturbed test function method \cite{perturbed}
to show that $\varphi$ is a viscosity solution of (\ref{eq:hamilton Jacobi}).
Theorem \ref{thm:Let--be} will follow by uniqueness of the solution
\cite{CEL}.

\subsection{Subsolution procedure}

Let $\psi\in\mathcal{C}^{1}\left(\mathbb{R}_{+}\times\mathbb{R}^{n}\right)$
be a test function such that $\varphi-\psi$ has a local strict maximum
at $\left(t^{0},x^{0}\right)$. We want to show that $\psi$ is a
subsolution of (\ref{eq:hamilton Jacobi}). If $\nabla_{x}\psi\left(t^{0},x^{0}\right)\in\mathrm{Sing}\left(M\right)^{c}$,
then we refer to \cite{BC12}.

Suppose now that $\nabla_{x}\psi\left(t^{0},x^{0}\right)\in\mathrm{Sing}\left(M\right)$.
Let $w\in\mathrm{Arg}\mu\left(\nabla_{x}\psi\left(t^{0},x^{0}\right)\right)\cap V$.
Then, $w\cdot\nabla_{x}\psi\left(t^{0},x^{0}\right)=\mu\left(\nabla_{x}\psi\left(t^{0},x^{0}\right)\right)$.
The uniform convergence of $\varphi^{\varepsilon}$ toward $\varphi$
ensures that the function $\left(t,x\right)\mapsto\varphi^{\varepsilon}\left(t,x,w\right)-\psi\left(t,x\right)$
has a local maximum at a point $\left(t^{\varepsilon},x^{\varepsilon}\right)$
satisfying $\left(t^{\varepsilon},x^{\varepsilon}\right)\to$$\left(t^{0},x^{0}\right)$,
as $\varepsilon\to0$. We then have:
\[
\partial_{t}\psi\left(t^{\varepsilon},x^{\varepsilon}\right)+w\cdot\nabla_{x}\psi\left(t^{\varepsilon},x^{\varepsilon}\right)=\partial_{t}\varphi^{\varepsilon}\left(t^{\varepsilon},x^{\varepsilon}\right)+w\cdot\nabla_{x}\varphi^{\varepsilon}\left(t^{\varepsilon},x^{\varepsilon}\right)=1-\int_{V} M\left(v'\right)e^{\frac{\varphi^{\varepsilon}\left(t^{\varepsilon},x^{\varepsilon},w\right)-\varphi{}^{\varepsilon}\left(t^{\varepsilon},x^{\varepsilon},v'\right)}{\varepsilon}}dv'\leq1.
\]

Passing to the limit $\varepsilon\to0$, we get \foreignlanguage{french}{$\partial_{t}\psi\left(t^{0},x^{0}\right)+\mu\left(\nabla_{x}\psi\left(t^{0},x^{0}\right)\right)\leq1$}.
We conclude that $\varphi$ is a viscosity subsolution of (\ref{eq:hamilton Jacobi}).

\subsection{Supersolution procedure}

Let $\psi\in\mathcal{C}^{1}\left(\mathbb{R}_{+}\times\mathbb{R}^{n}\right)$
be a test function such that $\varphi-\psi$ has a local strict minimum
at $\left(t^{0},x^{0}\right)$. We want to show that $\psi$ is a
supersolution of (\ref{eq:hamilton Jacobi}). If $\nabla_{x}\psi\left(t^{0},x^{0}\right)\in\mathrm{Sing}\left(M\right)^{c}$,
then we refer to \cite{BC12}.

Suppose now that $\nabla_{x}\psi\left(t^{0},x^{0}\right)\in\mathrm{Sing}\left(M\right)$.
Then, $\nabla_{x}\psi\left(t^{0},x^{0}\right)\neq0$ because $0\in\mathrm{Sing}\left(M\right)^{c}$.
We suppose without loss of generality that the minimum of $\varphi-\psi$
is global and that $\varphi\left(t^{0},x^{0}\right)-\psi\left(t^{0},x^{0}\right)=0$.
Let $\psi^{\varepsilon}:=\psi-C\left(t-t^{0}\right)^{2}+\varepsilon\eta$
with $C>0$ yet to be determined and
\[
\eta\left(v\right):=\mathrm{ln}\left(\mu\left(\nabla_{x}\psi\left(t^{0},x^{0}\right)\right)-v\cdot\nabla_{x}\psi\left(t^{0},x^{0}\right)\right).
\]
Then, $\eta$ is a continuous function on $D\left(\eta\right)=V\setminus\mathrm{Arg}\mu\left(\nabla_{x}\psi\left(t^{0},x^{0}\right)\right)$
and, for all $w\in\mathrm{Arg}\mu\left(\nabla_{x}\psi\left(t^{0},x^{0}\right)\right)\cap V$,
we have $\underset{v\to w}{\mathrm{lim}}\eta\left(v\right)=-\infty$.
Moreover, $\eta$ is bounded from below on all compact sets yielding the
uniform convergence $\psi^{\varepsilon}\to\psi$ on all compact sets
of $D\left(\eta\right)$. Finally, ${\int_{V}M\left(v'\right)e^{-\eta\left(v'\right)}dv'\leq1}$ since $\nabla_{x}\psi\left(t^{0},x^{0}\right)\in\mathrm{Sing}\left(M\right)$.

The function $\varphi-\left(\psi-C\left(t-t^{0}\right)^{2}\right)$
has a global strict minimum at $\left(t^{0},x^{0}\right)$. The first inequality
(\ref{apriori:Linfty}) ensures that the function $\varphi^{\varepsilon}-\psi^{\varepsilon}$
has a local minimum at a point $\left(t^{\varepsilon},x^{\varepsilon},v^{\varepsilon}\right)\in\mathbb{R}_{+}\times\mathbb{R}^{n}\times D\left(\eta\right)$.
As $V$ compact, we can extract a subsequence $\left(v^{\varepsilon}\right)_{\varepsilon}$, without relabelling, such that $v^{\varepsilon}\to v^{0}$, as $\varepsilon\to0$.

If $v^{0}\in V\setminus\mathrm{Arg}\mu\left(p\right)$, then there
exists a compact $A\subset D\left(\eta\right)$ such that $v^{0}\in A$
and the uniform convergence of $\psi^{\varepsilon}$ towards $\psi$
on $A$ guarantees that $\left(t^{\varepsilon},x^{\varepsilon}\right)\to\left(t^{0},x^{0}\right)$,
as $\varepsilon\to0$. We then get at point $\left(t^{\varepsilon},x^{\varepsilon},v^{\varepsilon}\right)$,\foreignlanguage{french}{
\begin{eqnarray*}
\partial_{t}\psi-2C\left(t^{\varepsilon}-t^{0}\right)+v^{\varepsilon}\cdot\nabla_{x}\psi=\partial_{t}\psi^{\varepsilon}+v^{\varepsilon}\cdot\nabla_{x}\psi^{\varepsilon}=\partial_{t}\varphi^{\varepsilon}+v^{\varepsilon}\cdot\nabla_{x}\varphi^{\varepsilon} & = & 1-\int_{V} M'e^{\frac{\varphi^{\varepsilon}-\varphi'^{\varepsilon}}{\varepsilon}}dv'\\
 & \geq & 1-\int_{V} M\left(v'\right)e^{\eta\left(v^{\varepsilon}\right)-\eta\left(v'\right)}dv'.
\end{eqnarray*}
}We take the limit $\varepsilon\to0$:\foreignlanguage{french}{
\begin{eqnarray*}
\partial_{t}\psi\left(t^{0},x^{0}\right)+v^{0}\cdot\nabla_{x}\psi\left(t^{0},x^{0}\right) & \geq & 1-e^{\eta\left(v^{0}\right)}\int_{V} M\left(v'\right)e^{-\eta\left(v'\right)}dv'\geq1-e^{\eta\left(v^{0}\right)}.
\end{eqnarray*}
}By construction, for all $v,v'\in D\left(\eta\right)$, we have $e^{\eta\left(v\right)}-e^{\eta\left(v'\right)}=\left(v'-v\right)\cdot\nabla_{x}\psi\left(t^{0},x^{0}\right)$
hence, for all $v\in D\left(\eta\right)$,we have $\partial_{t}\psi\left(t^{0},x^{0}\right)+v\cdot\nabla_{x}\psi\left(t^{0},x^{0}\right) \geq 1-e^{\eta\left(v\right)}$. Let $w\in V\cap\mathrm{Arg}\mu\left(\nabla_{x}\psi\left(t^{0},x^{0}\right)\right)$. Since $\mathrm{Arg}\mu\left(\nabla_{x}\psi\left(t^{0},x^{0}\right)\right)$ is a null-set, $V$ is dense in $\mathrm{Arg}\mu\left(\nabla_{x}\psi\left(t^{0},x^{0}\right)\right)$. Taking the limit $v\to w$, we get:
$\partial_{t}\psi\left(t^{0},x^{0}\right)+\mu\left(\nabla_{x}\psi\left(t^{0},x^{0}\right)\right) \geq 1.$

If $v^{0}\in V\cap\mathrm{Arg}\mu\left(p\right)$, we still have $\left(t^{\varepsilon},x^{\varepsilon}\right)\underset{\varepsilon\to0}{\longrightarrow}\left(t^{0},x^{0}\right)$
thanks to the following lemma:
\begin{lemma}
\label{lem:}For $C=4\left\Vert \varphi_{0}\right\Vert _{\infty}$, 
we have $\underset{\varepsilon\to0}{\mathrm{lim}}\,\varepsilon\eta\left(v^{\varepsilon}\right)=0$.\end{lemma}
\textbf{Proof of Lemma \ref{lem:}}
We have $\varphi^{\varepsilon}\left(t,x,v\right)-\varphi\left(t,x\right)\geq-2\left\Vert \varphi_{0}\right\Vert _{\infty} $ by (\ref{apriori:Linfty}) and
$\varphi\left(t,x\right)-\psi\left(t,x\right)\geq0$ hence
\begin{eqnarray*}
\varphi^{\varepsilon}\left(t,x,v\right)-\psi^{\varepsilon}\left(t,x,v\right) & \ge & -2\left\Vert \varphi_{0}\right\Vert _{\infty}+C\left(t-t^{0}\right)^{2}-\varepsilon\eta\left(v\right),\quad\forall\varepsilon>0.
\end{eqnarray*}
Moreover,
\begin{eqnarray*}
\varphi^{\varepsilon}\left(t^{0},x^{0},v\right)-\psi^{\varepsilon}\left(t^{0},x^{0},v\right) & = & \varphi^{\varepsilon}\left(t^{0},x^{0},v\right)-\varphi\left(t^{0},x^{0}\right)-\varepsilon\eta\left(v\right)\leq2\left\Vert \varphi_{0}\right\Vert _{\infty}-\varepsilon\eta\left(v\right).
\end{eqnarray*}
Since $C=4\left\Vert \varphi_{0}\right\Vert _{\infty}$, we have $\varphi^{\varepsilon}\left(t,x,v\right)-\psi^{\varepsilon}\left(t,x,v\right)>\varphi^{\varepsilon}\left(t^{0},x^{0},v\right)-\psi^{\varepsilon}\left(t^{0},x^{0},v\right)$
for all $t>t^{0}+1$ and, thus, the minimum of $\varphi^{\varepsilon}-\psi^{\varepsilon}$
cannot be attained for $t>t^{0}+1$ hence $t^{\varepsilon}\leq t^{0}+1$
for all $\varepsilon>0$. At point $\left(t^{\varepsilon},x^{\varepsilon},v^{\varepsilon}\right)$
we have:\foreignlanguage{french}{
\begin{eqnarray*}
\nabla_{v}\varphi^{\varepsilon}\left(t^{\varepsilon},x^{\varepsilon},v^{\varepsilon}\right) & = & \nabla_{v}\psi^{\varepsilon}\left(t^{\varepsilon},x^{\varepsilon},v^{\varepsilon}\right)=\varepsilon\nabla_{v}\eta\left(v^{\varepsilon}\right)=-\frac{\varepsilon\nabla_{x}\psi\left(t^{0},x^{0}\right)}{\mu\left(\nabla_{x}\psi\left(t^{0},x^{0}\right)\right)-v^{\varepsilon}\cdot\nabla_{x}\psi\left(t^{0},x^{0}\right)}.
\end{eqnarray*}
}The second estimation (\ref{apriori:Linfty}) yields $\left\Vert \nabla_{v}\varphi^{\varepsilon}\left(t^{\varepsilon},\cdot,\cdot\right)\right\Vert _{\infty}\leq t^{\varepsilon}\left\Vert \nabla_{x}\varphi_{0}\right\Vert _{\infty}\leq\left(t^{0}+1\right)\left\Vert \nabla_{x}\varphi_{0}\right\Vert _{\infty}$
hence
\begin{eqnarray*}
\frac{\varepsilon}{\left(t^{0}+1\right)\left\Vert \nabla_{x}\varphi_{0}\right\Vert _{\infty}}\left|\nabla_{x}\psi\left(t^{0},x^{0}\right)\right| & \leq & \mu\left(\nabla_{x}\psi\left(t^{0},x^{0}\right)\right)-v^{\varepsilon}\cdot\nabla_{x}\psi\left(t^{0},x^{0}\right),\\
\implies\varepsilon K\geq\varepsilon\eta\left(v^{\varepsilon}\right) & \geq & \varepsilon\mathrm{ln}\left(\frac{\varepsilon}{\left(t^{0}+1\right)\left\Vert \nabla_{x}\varphi_{0}\right\Vert _{\infty}}\left|\nabla_{x}\psi\left(t^{0},x^{0}\right)\right|\right),
\end{eqnarray*}
and $\varepsilon\eta\left(v^{\varepsilon}\right)\to0$ as $\varepsilon\to0$.\begin{flushright}$\square$\end{flushright}

Thanks to Lemma \ref{lem:}, the function $\left(t,x\right)\mapsto\psi^{\varepsilon}\left(t,x,v^{\varepsilon}\right)=\psi\left(t,x\right)-4\left\Vert \varphi_{0}\right\Vert _{\infty}\left(t-t^{0}\right)^{2}+\varepsilon\eta\left(v^{\varepsilon}\right)$
converges uniformly towards $\left(t,x\right)\mapsto\psi\left(t,x\right)-4\left\Vert \varphi_{0}\right\Vert _{\infty}\left(t-t^{0}\right)^{2}$
and has a local minimum at $\left(t^{\varepsilon},x^{\varepsilon}\right)$
satisfying $\left(t^{\varepsilon},x^{\varepsilon}\right)\to\left(t^{0},x^{0}\right)$,
as $\varepsilon\to0$. At point $\left(t^{\varepsilon},x^{\varepsilon},v^{\varepsilon}\right)$,
we have:
\[
\partial_{t}\psi^{\varepsilon}+v^{\varepsilon}\cdot\nabla_{x}\psi^{\varepsilon}=\partial_{t}\varphi^{\varepsilon}+v^{\varepsilon}\cdot\nabla_{x}\varphi^{\varepsilon}=1-\int_{V}M\left(v'\right)e^{\frac{\varphi^{\varepsilon}\left(t^{\varepsilon},x^{\varepsilon},v^{\varepsilon}\right)-\varphi^{\varepsilon}\left(t^{\varepsilon},x^{\varepsilon},v'\right)}{\varepsilon}}dv'.
\]
The minimal property of $\left(t^{\varepsilon},x^{\varepsilon},v^{\varepsilon}\right)$
implies at this point:
\begin{eqnarray*}
\partial_{t}\psi\left(t^{\varepsilon},x^{\varepsilon}\right)-8\left\Vert \varphi_{0}\right\Vert _{\infty}\left(t^{\varepsilon}-t^{0}\right)+v^{\varepsilon}\cdot\nabla_{x}\psi\left(t^{\varepsilon},x^{\varepsilon}\right)=\partial_{t}\psi^{\varepsilon}+v^{\varepsilon}\cdot\nabla_{x}\psi^{\varepsilon} & \geq & 1-\int_{V} M\left(v'\right)e^{\eta\left(v^{\varepsilon}\right)-\eta\left(v'\right)}dv'\\
 & \geq & 1-e^{\eta\left(v^{\varepsilon}\right)}.
\end{eqnarray*}
Passing to the limit $\varepsilon\to0$, we get $\partial_{t}\psi\left(t^{0},x^{0}\right)+\mu\left(\nabla_{x}\psi\left(t^{0},x^{0}\right)\right)\geq1$.
We conclude that $\varphi$ is a viscosity supersolution of (\ref{eq:hamilton Jacobi}).\begin{flushright}$\square$\end{flushright}

% \appendix

% \section{}
% \label{}

% The Acknowledgements are an un-numbered section
\section*{Acknowledgements}
This project has received funding from the European Research Council 
(ERC) under the European Union's Horizon 2020 research and innovation 
programme (grant agreement No 639638).

The author also wishes to thank Vincent Calvez and Julien Vovelle for their kind help.

\end{document}